\documentclass[12pt]{article}
\usepackage{amsmath}
\usepackage{amsthm}
\usepackage{graphicx}
\usepackage{verbatim}
\usepackage{epsfig}
\usepackage{hyperref}
\usepackage{float}
\usepackage{color}
\usepackage{fullpage}
\usepackage{amsfonts}
\usepackage{amscd}
\usepackage{amssymb}
\usepackage{graphics}
\usepackage{amsmath}
\usepackage[english]{babel}
\usepackage{bbm}
\usepackage{yfonts}
\usepackage{mathrsfs}
\usepackage{color}
\DeclareFontFamily{OML}{rsfs}{\skewchar\font'177}
\DeclareFontShape{OML}{rsfs}{m}{n}{ <5> <6> rsfs5 <7> <8> <9>
rsfs7 <10> <10.95> <12> <14.4> <17.28> <20.74> <24.88> rsfs10 }{}
\DeclareMathAlphabet{\mathfs}{OML}{rsfs}{m}{n}

\newcommand{\BE}{{\mathbb{E}}}

\newcommand{\BN}{{\mathbb{N}}}

\newcommand{\BP}{{\mathbb{P}}}

\newcommand{\BR}{{\mathbb{R}}}

\newcommand{\BZ}{{\mathbb{Z}}}

\newcommand{\CA}{{\mathcal{A}}}
\newcommand{\CB}{{\mathcal{B}}}
\newcommand{\CC}{{\mathcal{C}}}
\renewcommand{\CD}{{\mathcal{D}}}

\newcommand{\CO}{{\mathcal{O}}}

\newcommand{\bae}{\begin{equation}\begin{aligned}}
\newcommand{\eae}{\end{aligned}\end{equation}}

\newcommand{\psgeh}{\Psi_{\geq h}}
\newcommand{\psgh}{\Psi_{>h}}

\newtheorem{theorem}{Theorem}[section]
\newtheorem{proposition}[theorem]{Proposition}
\newtheorem{lemma}[theorem]{Lemma}

\newtheorem{corollary}[theorem]{Corollary}
\newtheorem{definition}[theorem]{Definition}

\begin{document}
\numberwithin{equation}{section} \numberwithin{figure}{section}
\title{Phase transition and uniqueness of levelset percolation}
\author{Erik Broman and Ronald Meester}
\maketitle

\begin{abstract}
The main purpose of this paper is to introduce and establish basic results 
of a natural extension of the classical Boolean 
percolation model (also known as the Gilbert disc model). 
We replace the balls of that model by a positive non-increasing attenuation function 
$l:(0,\infty) \to (0,\infty)$ to create the random field 
$\Psi(y)=\sum_{x\in \eta}l(|x-y|),$ where $\eta$ is a homogeneous 
Poisson process in $\BR^d.$ The field $\Psi$ is then a random potential field
with infinite range dependencies whenever the support of the function $l$ is unbounded.

In particular, we study the level sets $\Psi_{\geq h}(y)$ containing the 
points $y\in \BR^d$ such that $\Psi(y)\geq h.$
In the case where $l$ has unbounded support, we give, for any $d\geq 2,$
exact conditions on $l$ for $\Psi_{\geq h}(y)$ 
to have a percolative phase transition as a function of $h.$ 

We also prove that when $l$ is continuous then so is $\Psi$ almost surely. 
Moreover, in this case and for $d=2,$ we prove uniqueness of the 
infinite component of $\Psi_{\geq h}$ 
when such exists, and we also show that the so-called percolation function 
is continuous below the critical value $h_c$.

\end{abstract}

\section{Introduction}

In the classical Boolean continuum 
percolation model (see \cite{MR} for an overview), 
one considers 
a homogeneous Poisson process $\eta$ of rate $\lambda>0$ in $\BR^d$, 
and around each point $x\in \eta$ one places a ball $B(x,r)$ of radius $r.$ 
The main object of study is then 
\begin{equation} \label{eqn:Cdef}
\CC:=\bigcup_{x\in \eta}B(x,r),
\end{equation}
which is referred to as the {\em occupied} set.
It is well known (see \cite{MR}, Chapter 3) that there exists an 
$r_c=r_c(d)\in(0,\infty)$ such that 
\[
r_c:=\inf\{r:\BP(\CC(r) \textrm{ contains an unbounded component})>0\}.
\]
It is also well known that 
$\BP(\CC(r) \textrm{ contains an unbounded component})\in \{0,1\}.$
An immediate scaling argument shows that varying $\lambda$ is equivalent
to varying $r,$ and so one can fix $\lambda=1.$
This model was introduce by Gilbert in \cite{Gilbert} and further studied 
in \cite{Alexander}, \cite{BenSchramm}, \cite{MenshSid} and \cite{Roy}
(to name a few), while a dynamical version of this model was studied 
in \cite{ABGM}.

We consider a natural extension of this model. 
Let $\eta$ be a Poisson process with rate $\lambda$ in
$\BR^d,$ and let $x\in \eta$ denote a point in this process
(here we use the standard abuse of notation by writing 
$x\in \eta$ instead of $x\in {\rm supp}(\eta)$).
Furthermore, let $l:(0,\infty) \to (0,\infty)$ be a non-increasing function 
that we will call the {\em attenuation} function.
We then define the random field 
$\Psi=\Psi(l,\eta)$ at any point $y\in \BR^d$ by
\begin{equation} \label{eqn:Psidef}
\Psi(y):=\sum_{x\in \eta}l(|x-y|).
\end{equation}
In order for this to be well defined at every point, we let 
$l(0):=\lim_{r \to 0}l(r)$ (which can possibly be infinite). 
One can think of $\Psi$ as a random potential field where the 
contribution form a single point $x\in \eta$ is determined 
by the function $l.$

For any $0<h<\infty,$ we define 
\[
\Psi_{\geq h}:=\{y\in \BR^d: \Psi(y)\geq h\},
\]
which is simply the part of the random field $\Psi$  
which is at or above the level $h.$ Sometimes we also need 
\[
\Psi_{> h}:=\{y\in \BR^d: \Psi(y)> h\}.
\]

We note that if we consider our general model with $l(|x|)=I(|x|\leq r)$ 
(where $I$ is an indicator function), we have that 
$\CC$ and $\Psi_{\geq 1}$ have the same distribution, so the Boolean 
percolation model can be regarded 
as a special case of our more general model. 

When $l$ has unbounded support, adding or removing
a single point of $\eta$ will affect the field $\Psi$ at every point
of $\BR^d.$ Thus, our model does not have a so-called finite
energy condition 
which is the key 
to many standard proofs
in percolation theory. This is what makes studying $\Psi$
challenging (and in our opinion interesting). However, if we assume 
that $l$ has bounded support, a version of finite energy is recovered 
(see also the remark in Section \ref{sec:uniqueness} after the proof of 
our uniqueness result, Theorem \ref{thm:uniqueness}).

It is easy to see that varying $h$
and varying $\lambda$ is {\em not} equivalent. However, we will 
nevertheless restrict our attention to the case $\lambda=1$. 
In fact, there are many different sub-cases 
and generalizations that can be studied.
For instance: We can let $\lambda\in \BR,$ we can study $l$ having bounded
or unbounded support, we can let $l$ be a bounded or unbounded function,
let $l$ be continuous or discontinuous and we can study $\Psi_{\geq h}$
or $\Psi_{> h}$, to name a few possibilities.
While some results (Theorems \ref{thm:hcfinite} and \ref{thm:hcNT}) 
include all or most of the cases listed, others 
(Proposition \ref{prop:contfield} and Theorems \ref{thm:uniqueness}
and Theorem \ref{thm:thetacont}) require more specialized proofs.
The purpose of this paper is {\em not} to handle all different cases. 
Instead, we will focus on the extension of the classical Boolean percolation model
that we find to be the most natural
and interesting; when $l$ is continuous and with 
unbounded support.\\

We will now proceed to state our results, but first we 
have the following natural definition.
\begin{definition}
If $\Psi_{\geq h}$ ($\Psi_{> h}$) contains an unbounded connected component, 
we say that 
$\Psi$ {\em percolates} at (above) level $h,$ or simply that $\Psi_{\geq h}$
($\Psi_{> h}$) percolates.
\end{definition}
One would of course expect that percolation occurs either with probability 0 
or with probability 1, and indeed, our next result shows just that.
\begin{proposition} \label{prop:perc01}
We have that 
\[
\BP(\Psi_{\geq h} \textrm{ percolates})\in\{0,1\}
\textrm{ and } \BP(\Psi_{>h} \textrm{ percolates})\in\{0,1\}.
\]
\end{proposition}

\medskip\noindent
{\bf Proof.} 
This follows from a classical 
ergodicity argument. Indeed, the random field $\Psi$ is ergodic with 
respect to the 
group of translations of the space, see for instance the argument 
in \cite{MR}, Section 2.1, 
where it is formulated for the Boolean model and the random connection 
model, but the 
argument applies to our case as well. Since the event that 
$\Psi_{\geq h}$ percolates 
is invariant under translations, it must then have probability 0 or 1. 
\fbox{}\\

We now define
\[
h_c:=\sup\{h:\Psi_{\geq h} \textrm{ percolates with probability 1}\}.
\]
If we define $\tilde{h}_c$ as above, but with $\Psi_{\geq h}$
replaced by $\Psi_{>h},$ we see that $h_c=\tilde{h}_c.$ Indeed, 
if $h<h_c$ then $h\leq \tilde{h}_c$ so that $h_c\leq \tilde{h}_c.$
Also, if $h_c<\tilde{h}_c$ then we can find $h_c<g<\tilde{h}_c$
so that $\Psi_{>g}$ percolates while $\Psi_{\geq g}$ does not.
This is clearly impossible since $\Psi_{\geq g} \subset \Psi_{>g}.$

One of the main efforts of this paper is to establish conditions 
under which $h_c$ is nontrivial. 
As we will see, our results are qualitatively different depending on 
whether the attenuation function $l$ has bounded support or not. 
Our first main result is the following. 

\begin{theorem} \label{thm:hcfinite}
If the attenuation function $l$ satisfies $\int_{1}^\infty r^{d-1}l(r)dr<\infty,$
then $h_c<\infty.$ If instead $\int_{1}^\infty r^{d-1}l(r)dr=\infty,$
then almost surely $\Psi(y)=\infty$ for every $y\in \BR^d$, and so $h_c=\infty.$
\end{theorem}
\noindent
{\bf Remark:} The choice of the lower integral boundary 1 in $\int_{1}^\infty r^{d-1}l(r)dr$ 
is somewhat arbitrary, as replacing it with $\int_{c}^\infty r^{d-1}l(r)dr$ for any 
$0<c<\infty,$ would give the same result. 
Note also that if $l$ has bounded support, we must have that 
$\int_{1}^\infty r^{d-1}l(r)dr<\infty.$
\bigskip

Our next result concerns the positivity of $h_c.$ The proof is 
very straightforward and complements Theorem \ref{thm:hcfinite}. 
We define $r_l:=\sup\{r:l(r)>0\}$.

\begin{theorem} \label{thm:hcNT}
For $d=2,$ then $h_c>0$ iff $r_l> r_c.$
For $d\geq 3,$ $h_c>0$ if $r_l> r_c$ while
$h_c=0$ if $r_l< r_c.$
\end{theorem}
\noindent
{\bf Remark:} As is clear from the proof, the gap when $r_l=r_c$
for $d\geq 3,$ is simply due to the fact that for $d \geq 3$ it is unknown whether 
$\BP(\CC(r_c) \textrm{ contains an unbounded component})$ is 0 (as when $d=2)$ 
or 1.
 
\bigskip

We highlight our interest in the case when $l$ has unbounded support by formulating the following 
immediate corollary.
\begin{corollary}
If the attenuation function $l$ has unbounded support, 
then $0<h_c<\infty$ if and only if 
 $\int_{1}^\infty r^{d-1}l(r)dr<\infty.$
\end{corollary}
For the rest of this section, we will only consider functions $l$
satisfying $\int_1^\infty r^{d-1}l(r)dr<\infty$.

We remark that mere a.s.\ finiteness of the field $\Psi$ does not guarantee that 
there is a nontrivial phase transition. Indeed, one can construct an example 
(following the examples of \cite{MR} Chapter 7.7)
of a stationary process together with a suitable attenuation function so that the 
ensuing field is finite a.s., while $h_c=\infty.$ 

Next, we turn to the everywhere continuity of the field $\Psi$. 
Of course, if $l$ is not continuous, everywhere continuity of $\Psi$ cannot hold. 
However, if $l$ is continuous then everywhere 
continuity of the field $\Psi$ would be expected. We note that in the case of 
$l$ being unbounded (i.e. $\lim_{r \to 0}l(r)=\infty$), we simply define
$\Psi(y)=\infty$ for every $y\in \eta.$

\begin{proposition} \label{prop:contfield}
If $l$ is continuous, then
the random field $\Psi$ is a.s.\ everywhere continuous.
\end{proposition}

Proposition \ref{prop:contfield}
will be of use when proving our next result concerning uniqueness of the 
unbounded component.

\begin{theorem}
\label{thm:uniqueness}
Let $h$ be such that $\Psi_{\geq h}$ ($\Psi_{>h}$) contains an 
unbounded component.
If $l$ is continuous and with unbounded support, then for $d=2$,
there is a unique such unbounded component.
\end{theorem}
\noindent
{\bf Remarks:} 
We will prove this theorem first for $\Psi_{>h}$ and then infer it for 
$\Psi_{\geq h}$, see also the discussion before the proof of the theorem.

There are of course a number of possible generalizations of 
this statement, and perhaps the most interesting/natural would be to investigate
it for $d\geq 3.$ We discuss this in some detail after the proof of 
Theorem \ref{thm:uniqueness}. \\

Let $\CC_{o,\geq}(h)$ ($\CC_{o,>}(h)$) be the connected component of 
$\Psi_{\geq h}$ ($\Psi_{> h}$) that contains the origin $o.$ Define 
the percolation function
\[
\theta_{\geq}(h):=\BP(\CC_{o,\geq}(h) \textrm{ is unbounded}),
\]
and similarly define $\theta_>(h)$. Our last result is the following.

\begin{theorem} \label{thm:thetacont}
The functions $\theta_{\geq}(h)$ and $\theta_{>}(h)$ are equal and 
continuous for $h<h_c$.
\end{theorem}

\medskip

The rest of the paper is organized as follows. 
In Section \ref{sec:phasetransition}
we prove Theorems \ref{thm:hcfinite} and \ref{thm:hcNT}. The continuity of 
$\Psi$ (Proposition \ref{prop:contfield})
is proved in Section \ref{sec:cont} and this results is then used in 
Section \ref{sec:uniqueness} in order to prove Theorem \ref{thm:uniqueness},
which in turn will allow us to prove Theorem \ref{thm:thetacont}.

\section{Conditions for the non-triviality of $h_c$} \label{sec:phasetransition}
We start by proving Theorem \ref{thm:hcNT} as this is easily done. 

\medskip
{\bf Proof of Theorem \ref{thm:hcNT}.}
Recall the constructions \eqref{eqn:Cdef} and \eqref{eqn:Psidef}, and observe 
that by using the same realization of $\eta,$ we can couple $\Psi$ and $\CC$
on the same probability space in the natural way. By using this coupling, we see that 
for any $h>0,$ $\Psi_{\geq h}\subset \CC(r_l).$ In the case $d=2,$
it is known (see \cite{MR}, Theorem 4.5) that $\CC(r_c)$ does not 
percolate, showing that $h_c=0$ when $d=2$ and $r_l\leq r_c.$
For $d\geq 3,$ the statement follows by 
observing that $\CC(r)$ does not percolate for $r<r_c$ by definition 
of $r_c.$

Assume instead that $r_l> r_c.$ Let $r_c<r<r_l,$ and  $h=h(r)$
be any $h$ such that 
\[
B(0,r) \subset \{y:l(|y|)\geq h\}.
\]
With this choice of $h,$ we see that any point $y$ within radius $r$ from some 
$x\in \eta$ will also belong to $\Psi_{\geq h}$ so that
\[
\CC(r)\subset \Psi_{\geq h}.
\]
Since $r>r_c,$ $\CC(r)$ a.s.\ contains an unbounded component and hence so does 
$\Psi_{\geq h}.$
\fbox{}\\

The proof of Theorem \ref{thm:hcfinite} is much more involved, and will require
a number of preliminary lemmas to be established first. In order to see what the 
purpose of these will be, we start by giving an outline of the 
strategy of our proof along with introducing some of the relevant notation.
Let $\alpha \BZ^d$ denote the lattice with 
spacing $\alpha>0.$ For any $z\in \alpha \BZ^d,$ let $B(z,\alpha)$
denote the closed box of side length $\alpha$ centered at $z,$ and define
$\CB_\alpha:=\{B(z,\alpha):z\in \alpha \BZ^d\}.$
For convenience, we assume from now on that $\alpha<1.$ 
 
\medskip\noindent
{\bf Claim:} There exists an $\epsilon>0$ such that if for 
any $0<\alpha<1$ and every $k$ and 
collection of distinct cubes $B_1,\ldots,B_k\in \CB_\alpha,$ we have that
\begin{equation} \label{eqn:epsest}
\BP(\sup_{y\in B_1}\Psi(y)\geq h, \ldots, \sup_{y\in B_k}\Psi(y)\geq h)
\leq \epsilon^k,
\end{equation}
then $\Psi_{\geq h}$ does a.s.\ not contain an unbounded component. \\

This claim can be proved using standard percolation arguments as follows. 
Let $B^o\in \CB_\alpha$ be the cube containing the origin $o$ and let 
$\CO$ denote the event that $B^o$ intersects an
unbounded component of $\Psi_{\geq h}$. If $\CO$ occurs, then
for any $k,$ there must exist a sequence
$B_1,B_2,\ldots B_k\in \CB_\alpha$ such that $B_1=B^o,$
$B_i \neq B_j$ for every $i\neq j,$  $B_i\cap B_{i+1}\neq \emptyset$
for every $i=1,\ldots,k-1$ and with the property that 
$\sup_{y\in B_i}\Psi(y)\geq h$ for every $i=1,\ldots,k.$
We note that the number of such paths must be bounded
by $3^{dk}$, as any box has fewer than $3^d$ 'neighbors'. 
Thus, from \eqref{eqn:epsest} we get that 
$\BP(\CO)\leq 3^{dk}\epsilon^k,$ and since this holds for
arbitrary $k$ this proves the claim by taking $\epsilon<3^{-d}$.

One issue when proving \eqref{eqn:epsest} is that we want
to consider the supremum of the field within the boxes $B_1,\ldots,B_k.$
However, this is fairly easily dealt with by introducing an auxiliary 
field $\tilde{\Psi}$ with the property that for any $B\in \CB_\alpha$
$\tilde{\Psi}(y_c(B))\geq \sup_{y\in B} \Psi(y)$ where $y_c(B)$
denotes the center of $B$ (see further \eqref{eqn:Psisup}). This allows us
to consider $k$ fixed points of the new field $\tilde{\Psi}$ rather than 
the supremums involved in \eqref{eqn:epsest}.

One of the main problems in proving \eqref{eqn:epsest} is the 
long range dependencies involved whenever $l$ has unbounded support 
(as discussed in the introduction). The strategy to resolve this issue
is based on the simple observation that 
\begin{equation} \label{eqn:supsumbound}
\left\{\sup_{y\in B_1}\Psi(y)\geq h, \ldots, \sup_{y\in B_k}\Psi(y)\geq h\right\}
\subset \left\{ \sum_{l=1}^k \tilde{\Psi}(y_c(B_l))\geq kh \right\}.
\end{equation}
The event on the right hand side of \eqref{eqn:supsumbound} can be analyzed 
using a version of Campbell's theorem (see e.g.\ \cite{Kingman} p. 57-57). An obvious
problem with this is that if $l$ is unbounded and if a single point of 
$\eta$ falls in $\bigcup_{l=1}^k B_l$, then the sum in \eqref{eqn:supsumbound} is 
infinite. However, by letting $\alpha$ above be very small, we can make sure that 
with very high probability, ``most'' of the boxes $B_1,\ldots,B_k$ will 
not contain any points of $\eta$ (and in fact there will not even be a point
in a certain neighborhood of the box). 
We then use a more sophisticated version of 
\eqref{eqn:supsumbound} (i.e. \eqref{eqn:prelest}) where we condition 
on which of the boxes $B_1,\ldots,B_k$ have a point of $\eta$ in their neighborhood,
and then sum only over the boxes whose neighborhoods are vacant of points.
This in turn introduces another problem, namely that we now have to deal 
with a Poisson process conditioned on the presence and absence of points
of $\eta$ in the neighborhoods of the boxes $B_1,\ldots,B_k.$ 
In particular, we have to control the
damage from knowing the presence of such points. This is the purpose of 
Lemmas \ref{lemma:elemPoisdom} and \ref{lemma:elemPoisdom2}, which will tell
us that our knowledge is not worse than having no information at all plus
adding a few extra points to the process. Later, Lemma 
\ref{lemma:bndedfield} will enable us to control the effect of this addition 
of extra points.

\bigskip

We now start presenting the rigorous proofs. Our first lemma 
is elementary, and the result is presumably folklore. However,
we give a proof for sake of completeness.
\begin{lemma} \label{lemma:elemPoisdom}
Let $X$ be a Poisson distributed random variable with parameter $\lambda.$ 
We have that for any $k\geq 0,$
\begin{equation}\label{eqn:Poiconddom}
\BP(X\geq k | X\geq 1)\leq \BP(X\geq k-1).
\end{equation}
\end{lemma}

\noindent
{\bf Proof.}
We claim that for any $X^n\sim$ Bin($n,p$) where $np=\lambda,$ and 
any $0\leq k \leq n,$ we have that 
\begin{equation} \label{eqn:Binconddom}
\BP(X^n\geq k| X^n \geq 1)\leq \BP(X^n\geq k-1).
\end{equation}
We observe that from \eqref{eqn:Binconddom} we get that 
\[
\BP(X\geq k | X\geq 1)
=\lim_{n \to \infty}\BP(X^n\geq k| X^n \geq 1)
\leq \lim_{n \to \infty}\BP(X^n\geq k-1)= \BP(X \geq k-1).
\]
This establishes \eqref{eqn:Poiconddom}, and so we need to prove 
\eqref{eqn:Binconddom}.

We will prove \eqref{eqn:Binconddom} through induction, and we start by observing
that it trivially holds for $n=1$ and $k=0,1.$ Assume therefore that 
\eqref{eqn:Binconddom} holds for $n$ and any $k=0,\ldots,n.$
We will write $X^n=X_1+\ldots+X_{n}$ where $\left(X_i\right)_{i\geq 1}$ is an 
i.i.d. sequence with $\BP(X_i=1)=1-\BP(X_i=0)=p.$ Of course, 
\eqref{eqn:Binconddom} trivially holds for $n+1$ and $k=0.$ Furthermore, 
we have that for any $k=1,\ldots,n,$
\begin{eqnarray*}
\lefteqn{\BP(X^{n+1} \geq k | X^{n+1}\geq 1)}\\
& & =\BP(X^{n+1} \geq k | X^{n+1}\geq 1, X_{n+1}=1)\BP(X_{n+1}=1 |X^{n+1}\geq 1)\\
& & \ \ \ \ +\BP(X^{n+1} \geq k | X^{n+1}\geq 1, X_{n+1}=0)\BP(X_{n+1}=0 |X^{n+1}\geq 1)\\
& & =\BP(X^{n} \geq k-1)\BP(X_{n+1}=1 |X^{n+1}\geq 1)
 +\BP(X^{n} \geq k | X^{n}\geq 1)\BP(X_{n+1}=0 |X^{n+1}\geq 1)\\
& & \leq \BP(X^{n} \geq k-1),
\end{eqnarray*}
where we use the induction hypothesis that $\BP(X^{n} \geq k | X^{n}\geq 1)
\leq \BP(X^{n} \geq k-1)$ in the last inequality. Finally, 
\[
\BP(X^{n+1} =n+1 | X^{n+1}\geq 1)
=\BP(X^{n} = n)\BP(X_{n+1}=1 |X^{n+1}\geq 1)
\leq \BP(X^{n} = n),
\]
and this establishes \eqref{eqn:Binconddom} for $n+1$ and any $k=0,\ldots,n+1.$
\fbox{}\\

Let $A_1,A_2,\ldots,A_n$ be subsets of $\BR^d, $ and let $C_1,\ldots,C_m$
be a partition of $\cup_{i=1}^n A_i$ such that for any $i,$ 
\begin{equation} \label{eqn:Aipart}
A_i=\cup_{k=1}^l C_{i_k},
\end{equation}
for some collection $C_{i_1},\ldots,C_{i_l}.$ Let $\eta_A$ be a homogeneous 
Poisson process of rate $\lambda>0$ on $\cup_{i=1}^n A_i$ conditioned 
on the event $\cap_{i=1}^n\{\eta(A_i)\geq 1\},$ and 
let $\eta'_A$ be a homogeneous (unconditioned) Poisson process of 
rate $\lambda>0$ on $\cup_{i=1}^n A_i.$ Furthermore, let 
$\xi_A$ be a point process on $\cup_{i=1}^n A_i$ consisting of
exactly one point in each of the sets $C_1,\ldots,C_m$ such that 
the position of the point in $C_i$ is uniformly distributed within 
the set, and so that this position is independent between sets.

Our next step is to use Lemma \ref{lemma:elemPoisdom} to prove 
a result relating the conditioned Poisson process $\eta_A$ to 
the sum $\eta'_A+\xi_A,$ where $\eta'_A$ and $\xi_A$ are independent.
For two point processes $\eta_1,\eta_2$ in $\BR^d$, we write $\eta_1\preceq \eta_2$ if there exists a coupling 
of $\eta_1,\eta_2$ so that $\BP(\eta_1\subset \eta_2)=1.$
\begin{lemma} \label{lemma:elemPoisdom2}
Let $\eta_A,\eta'_A$ and $\xi_A$ be as above, and let 
$\eta'_A$ and $\xi_A$ be independent. We have that 
\[
\eta_A \preceq \eta'_A+\xi_A.
\]
\end{lemma}

Informally, Lemma \ref{lemma:elemPoisdom2} tells us that if 
we consider a homogeneous Poisson process conditioned on the presence
of points in $A_1,\ldots,A_k$, it is not worse than taking an unconditioned
process and adding single points to all the sets $C_1,\ldots,C_m$ (which are 
used as the building blocks for the sets $A_1,\ldots,A_k$).

\medskip

\noindent
{\bf Proof of Lemma \ref{lemma:elemPoisdom2}.}
As usual, let $\eta$ be a homogeneous Poisson process on $\BR^d.$
Let $J=(j_1,\ldots,j_m)\in\{0,1\}^m$ and define
\[
\CC_J=\bigcap_{l:j_l=1}\{\eta(C_l)\geq 1\}\bigcap_{l:j_l=0}\{\eta(C_l)=0\}.
\]
We note that either $\CC_J \subset \cap_{i=1}^n\{\eta(A_i)\geq 1\},$ or 
$\CC_J \cap_{i=1}^n\{\eta(A_i)\geq 1\}=\emptyset,$ which follows from 
\eqref{eqn:Aipart}. Indeed, if for any $i$ and $J\in \{0,1\}^m,$ 
all of the sets $C_{i_k}$ in 
\eqref{eqn:Aipart} have $\eta(C_{i_k})=0$, then 
$\CC_J \cap_{i=1}^n\{\eta(A_i)\geq 1\}=\emptyset.$ On the other hand, if for 
every $i,$ there exists some set $C_{i_k}$ in \eqref{eqn:Aipart} such that 
$\eta(C_{i_k})=1$ this implies that $\{\eta(A_i)\geq 1\}$ 
occurs for every $i.$

Using this, we have that for any $(k_1,\ldots,k_m)\in \BN^m,$ 
\begin{eqnarray} \label{eqn:CJineq}
\lefteqn{\BP(\eta(C_1)\geq k_1,\ldots, \eta(C_m)\geq k_m 
|\cap_{i=1}^n\{\eta(A_i)\geq 1\})} \\
& & =\sum_{J\in \{0,1\}^m}\BP(\eta(C_1)\geq k_1,\ldots, \eta(C_m)\geq k_m 
|\CC_J)\BP(\CC_J | \cap_{i=1}^n\{\eta(A_i)\geq 1\}), \nonumber
\end{eqnarray}
since $\BP(\CC_J | \cap_{i=1}^n\{\eta(A_i)\geq 1\})=0$ if 
$\CC_J \not \subset \cap_{i=1}^n\{\eta(A_i)\geq 1\}.$
Furthermore, for any $J\in \{0,1\}^m$
\begin{equation} \label{eqn:CJineq2}
\BP(\eta(C_1)\geq k_1,\ldots, \eta(C_m)\geq k_m |\CC_J)
\leq \BP(\eta(C_1)\geq k_1-1,\ldots, \eta(C_m)\geq k_m-1),
\end{equation}
by using Lemma \ref{lemma:elemPoisdom} and a trivial bound.

Combining \eqref{eqn:CJineq} and \eqref{eqn:CJineq2} yields
\begin{eqnarray*}
\lefteqn{\BP(\eta_A(C_1)\geq k_1,\ldots, \eta_A(C_m)\geq k_m )}\\
& & =\BP(\eta(C_1)\geq k_1,\ldots, \eta(C_m)\geq k_m 
|\cap_{i=1}^n\{\eta(A_i)\geq 1\})\\
& & \leq 
\BP(\eta(C_1)\geq k_1-1,\ldots, \eta(C_m)\geq k_m-1)\\
& & =\BP((\eta'_A+\xi_A)(C_1)\geq k_1,\ldots, (\eta'_A+\xi_A)(C_m)\geq k_m).
\end{eqnarray*}
The statement follows by the elementary property of a Poisson process,
that conditioned on a certain number of points falling within a fix set
D, these points are independently and uniformly distributed within that set.
\fbox{}\\

We now turn to the issue of taking the
supremum of the field over a box. Therefore, let $0<\alpha<1,$ and
define the auxiliary attenuation function $\tilde{l}_\alpha$ by 
\[
\tilde{l}_\alpha(r)=
\left\{
\begin{array}{ll}
l(0) & \textrm{if } r\leq \alpha \sqrt{d}/2 \\
l(r-\alpha\sqrt{d}/2) & \textrm{if } r\geq \alpha\sqrt{d}/2,
\end{array}
\right.
\]
for every $r\geq 0.$ If $y_{c}(B)$ denotes the center of the box $B\in \CB_\alpha,$
we note that for any $y\in B$ and $x\in \BR^d,$ 
\[
\tilde{l}_\alpha(|x-y_c(B)|)
\geq \tilde{l}_\alpha(|x-y|+|y-y_c(B)|)
\geq \tilde{l}_\alpha(|x-y|+\alpha\sqrt{d}/2)
= l(|x-y|).
\]
Therefore, if we let $\tilde{\Psi}$ be the field we get by using $\tilde{l}$
in place of $l$ in \eqref{eqn:Psidef}, we get that 
\begin{equation}\label{eqn:Psisup}
\tilde{\Psi}(y_c(B))=\sum_{x\in \eta}\tilde{l}_\alpha(|x-y_c(B)|)
\geq \sup_{y\in B}\sum_{x\in \eta} l(|x-y|)=\sup_{y\in B} \Psi(y).
\end{equation}

Our next lemma will be a central ingredient of the proof of 
Theorem \ref{thm:hcfinite}. It will deal with the effect to the field 
$\tilde{\Psi}$ of adding extra points to $\eta$. To that end, 
let $A_o$ be the box of side length $\alpha(4\lceil\sqrt{d}\rceil+1)$ 
centered around the origin $o$. For any box $B\in \CB_\alpha$ with 
$B\cap A_o=\emptyset,$
place a point $x_B$ in $B$ at the closest distance to the origin, and let 
$\xi$ denote the corresponding (deterministic) point set. Let
\[
\tilde{\Psi}_{A_o}(y):=\sum_{x \in \xi}\tilde{l}_\alpha(|x-y|),
\]
be the corresponding deterministic field. 

\begin{lemma} \label{lemma:bndedfield}
There exists a constant $C<\infty$ depending on $d$ but not on $\alpha$ and
such that for every $0<\alpha<1,$
\[
\tilde{\Psi}_{A_o}(o)\leq \frac{C}{\alpha^d} I_\alpha,
\]
where 
\[
I_\alpha=\int_{\alpha/2}^\infty r^{d-1} l(r)dr<\infty.
\]
\end{lemma}
\noindent
{\bf Proof.}
Consider some $B\in \CB_\alpha$ such that $B\cap A_o=\emptyset.$
We have that 
\[
\tilde{l}_\alpha(|x_B|)\leq \frac{1}{Vol(B)}\int_{B}\tilde{l}_\alpha(|x|-{\rm diam}(B))dx
=\frac{1}{\alpha^d}\int_{B}\tilde{l}_\alpha(|x|-\alpha \sqrt{d})dx.
\]
Therefore, 
\begin{eqnarray*}
\lefteqn{\tilde{\Psi}_{A_o}(o)
\leq \frac{1}{\alpha^d}\int_{\BR^d\setminus A_o}\tilde{l}_\alpha(|x|-\alpha \sqrt{d})dx}\\
& & \leq \frac{C}{\alpha^d}\int_{\alpha(2\lceil\sqrt{d}\rceil+1/2)}^\infty
r^{d-1}l(r-2\alpha \sqrt{d})dr \\
& & \leq \frac{C}{\alpha^d}\int_{\alpha/2}^\infty (r+\alpha(2\lceil\sqrt{d}\rceil))^{d-1}
l(x) dr \nonumber \\
& & \leq \frac{C}{\alpha^d}\int_{\alpha/2}^\infty (r+r(4\sqrt{d}+2))^{d-1}
l(r) dr \nonumber \\
& & = \frac{C}{\alpha^d}\int_{\alpha/2}^\infty r^{d-1} l(r)dr
\end{eqnarray*}
where the constant $C=C(d)<\infty$
is allowed to vary in the steps of the calculations. 
Finally, the fact that $I_\alpha<\infty,$ follows easily from the fact that 
$\int_{1}^\infty r^{d-1} l(r)dr<\infty.$
\fbox{}\\

We have now established all necessary tools in order to prove Theorem 
\ref{thm:hcfinite}. However, since the proofs of the two statements of 
Theorem \ref{thm:hcfinite} are very different, we start by proving the 
first one as a separate result.

\begin{theorem}\label{thm:hcfinite_aux}
If $\int_1^\infty l^{d-1}l(r)<\infty,$ then $h_c<\infty.$
\end{theorem}
\noindent
{\bf Proof.}
We shall prove that for any $\epsilon>0,$ \eqref{eqn:epsest} holds
for $\alpha$ small enough and $h$ large enough. This will prove our
result as explained just below \eqref{eqn:epsest}.

For any $B\in \CB_\alpha,$ let $A_\alpha(B)$ be the box concentric 
to $B$ and with side length $\alpha(4\lceil \sqrt{d}\rceil+1)$.
Let $E(B)$ be the event that $\eta(A_\alpha(B))=0,$ and observe that 
if $c=\BP(E(B)),$
we have that $c=c(\alpha)\to 1$ as $\alpha \to 0$. We say that the box 
$B$ is {\em good} if the event $E(B)$ occurs. Goodness of the boxes 
$B\in \CB_\alpha$ naturally induces a percolation model on $\BZ^d$ with 
a finite range dependency. Since the marginal probability $c(\alpha)$ of being 
good can be made to be arbitrarily close to 1 by taking 
$\alpha$ small enough, we can use Theorem B26 of \cite{SIS}
to dominate an i.i.d.\ product measure with density $p=p(\alpha)$ on 
the boxes $B\in \CB_\alpha.$
Furthermore, by the same theorem, we can take $p(\alpha)\to 1$ as $\alpha \to 0.$

Fix $k$ and a collection $B_1,B_2,\ldots, B_k$ as in \eqref{eqn:epsest}. 
For any $B_i,$ let $A_i=A_\alpha(B_i)$, and let $\Gamma_i:=I(B_i)$ 
where $I$ denotes an indicator function.
If we take $\Gamma=\sum_{i=1}^k \Gamma_i,$ then
by the above domination of a product measure of density $p,$ 
we see that $\Gamma$ is 
stochastically larger than $\Gamma'\sim$Bin$(p,k).$ Furthermore, 
we have that 
\begin{eqnarray*}
\lefteqn{\BP\left(\Gamma'\leq \frac{k}{2}\right)}\\
& & =\BP\left(e^{\log(1-p)\Gamma'}\geq e^{\log(1-p)k/2}\right)
\leq e^{-\log(1-p)k/2}\BE\left[e^{\log(1-p)\Gamma'}\right]\\
& & =e^{-\log(1-p)k/2}\left(pe^{\log(1-p)}+1-p\right)^k
\leq 2^k e^{\log(1-p)k/2}=e^{-d(\alpha)k},
\end{eqnarray*}
where we can take $d(\alpha)\to \infty$ as $\alpha \to 0,$ by 
taking $p(\alpha)\to 1.$
If we define $G_k$ to be the event that at least $k/2$ of the boxes 
$B_1,B_2,\ldots, B_k$
are good, we thus have that $\BP(G_k)\geq 1-e^{-d(\alpha)k}.$

Let $J=J(\eta)\in\{0,1\}^k$ be such that $J_j=1$ iff $B_j$ is good, and identify 
$J$ with the corresponding subset of $\{1,\ldots,k\}.$ Thus we write $j\in J$
iff $B_j$ is good. For any fixed
$J\in \{0,1\}^k,$ we let $\CD_J$ denote the event 
\[
\bigcap_{j\in J} E_j \bigcap_{j\in J^c} E_j^c
\]
so that $\CD_J$ is the event that each set $A_j$ such that $j\in J$
is vacant of points, while each set $A_j$ such that $j\in J^c$
contains at least one point of $\eta$.
We then have that 
\begin{eqnarray} \label{eqn:prelest}
\lefteqn{\BP(\sup_{y\in B_1}\Psi(y)\geq h, \ldots, \sup_{y\in B_k}\Psi(y)\geq h)}\\
& & \leq \sum_{J\in \{0,1\}^k}
\BP(\tilde{\Psi}(y_c(B_1))\geq h,\ldots,\tilde{\Psi}(y_c(B_k))\geq h|\CD_J)\BP(\CD_J)\nonumber \\
& & \leq \sum_{|J|\geq k/2}
\BP(\tilde{\Psi}(y_c(B_1))\geq h,\ldots,\tilde{\Psi}(y_c(B_k))\geq h|\CD_J)
\BP(\CD_J)+e^{-d(\alpha)k}, \nonumber 
\end{eqnarray}
by using \eqref{eqn:Psisup} in the first inequality and that 
$\BP(|J|<k/2)=\BP(G_k^c)$ in the last.

For any $J\in \{0,1\}^k,$ we let $\eta_J$ be a Poisson process on 
$\BR^d$ of rate 1, 
conditioned on the event $\CD_J.$ Furthermore, for $j\in J^c,$ let 
\[
D_j:=A_j \setminus \bigcup_{i\in J} A_i.
\]
We see that $\eta_J$ can be expressed as the sum of $\tilde{\eta}$ 
and $\eta_D$ where $\tilde{\eta}$ is a Poisson process on 
$\BR^d \setminus \bigcup_{n=1}^k A_n$ of unit rate, while $\eta_D$ is 
a Poisson process on $\bigcup_{j\in J^c} D_j$ conditioned on the 
event $\bigcap_{j\in J^c} \{\eta(D_j)\geq 1\}.$ We let
$C_1,\ldots,C_m\in \CB_\alpha$ be a partition (up to sets of measure zero) of $\bigcup_{j\in J^c} D_j$,
and use Lemma \ref{lemma:elemPoisdom2} to see that 
$\eta_D\preceq \eta_D'+\xi_D$. Here of course, $\eta'_D$ is Poisson process
on $\bigcup_{j\in J^c} D_j$ while $\xi_D$ is a point process consisting of 
one point added uniformly and independently to every box 
$C_i$ for $i=1,\ldots,m.$ As in Lemma \ref{lemma:elemPoisdom2} 
$\eta'_D$ and $\xi_D$ are independent.

We conclude that 
$\eta_J=\eta_D+\tilde{\eta} \preceq \eta'_D+\xi_D+\tilde{\eta}$,
and observe that $\eta'_D+\tilde{\eta}$ is just a homogeneous Poisson
process on $\BR^d\setminus A_J$ where $A_J:=\bigcup_{j\in J} A_j.$
Writing $\eta_{A_J^c}$ for this sum, we have that 
$\eta_J\preceq \xi_D+\eta_{A_J^c}$ and by
first using this, and then Markov's inequality, we get that
\begin{eqnarray} \label{eqn:Markov}
\lefteqn{\BP(\tilde{\Psi}(y_c(B_1))\geq h,\ldots,\tilde{\Psi}(y_c(B_k))\geq h
|\CD_J)}\\
& & \leq \BP_{\xi_D+\eta_{A_J^c}}
(\tilde{\Psi}(y_c(B_1))\geq h,\ldots,\tilde{\Psi}(y_c(B_k))\geq h) \nonumber \\
& & \leq \BP_{\xi_D+\eta_{A_J^c}}
\left(s\sum_{j\in J}\tilde{\Psi}(y_c(B_j))\geq h|J|s\right)
\leq e^{-h|J|s}\BE_{\xi_D+\eta_{A_J^c}}
\left[e^{s\sum_{j\in J}\tilde{\Psi}(y_c(B_j))}\right], \nonumber
\end{eqnarray}
where $\BP_{\xi_D+\eta_{A_J^c}}$ is the probability measure 
corresponding to $\xi_D+\eta_{A_J^c},$ and where 
$\BE_{\xi_D+\eta_{A_J^c}}$ denotes expectation with 
respect to this probability measure. We have
\begin{eqnarray*}
\lefteqn{\tilde{\Psi}(y_c(B_j))=\sum_{x\in\xi_D+\eta_{A_J^c}}
\tilde{l}_\alpha(|y_{c}(B_j)-x|)}\\
& & =\sum_{x\in \eta_{A_J^c}}\tilde{l}_\alpha(|y_{c}(B_j)-x|)
+\sum_{x\in \xi_D}\tilde{l}_\alpha(|y_{c}(B_j)-x|) \\
& & =\tilde{\Psi}_{\eta_{A_J^c}}(y_c(B_j))
+\tilde{\Psi}_{\xi_D}(y_c(B_j)),
\end{eqnarray*}
using obvious notation. Thus, using independence, we have that 
\begin{equation} \label{eqn:etaxi}
\BE_{\xi_D+\eta_{A_J^c}}
\left[e^{s\sum_{j\in J}\tilde{\Psi}(y_c(B_j))}\right]
=\BE_{\eta_{A_J^c}}
\left[e^{s\sum_{j\in J}\tilde{\Psi}_{\eta_{A_J^c}}(y_c(B_j))}\right]
\BE_{\xi_D}\left[e^{s\sum_{j\in J}\tilde{\Psi}_{\xi_D}(y_c(B_j))}\right].
\end{equation}
Consider the function $g_J(y):=\sum_{j\in J} \tilde{l}_\alpha(|y_{c}(B_j)-y|),$ 
so that
\[
\sum_{x\in \eta_{A_J^c}} g_J(x)
=\sum_{x\in \eta_{A_J^c}}\sum_{j\in J} \tilde{l}_\alpha(|y_{c}(B_j)-x|)
=\sum_{j\in J}\sum_{x\in \eta_{A_J^c}}\tilde{l}_\alpha(|y_{c}(B_j)-x|)
=\sum_{j\in J}\tilde{\Psi}_{\eta_{A_J^c}}(y_{c}(B_j)),
\]
and similarly 
\[
\sum_{x\in \xi_D} g_J(x)=\sum_{j\in J}\tilde{\Psi}_{\xi_D}(y_{c}(B_j)).
\]
By Campbell's theorem (see \cite{Kingman} p 57-58, \cite{D-VJ} Sections 2.4 and 9.4)
we get that
\begin{equation} \label{eqn:cambeta}
\BE_{\eta_{A_J^c}}
\left[e^{s\sum_{j\in J} \tilde{\Psi}_{\eta'_D+\tilde{\eta}}(y_{c}(B_j))}\right]
=\exp\left({\int_{\BR^d\setminus A_J}e^{sg_J(x)}-1 dx}\right).
\end{equation}
We proceed by bounding the right hand side of this expression and 
start by noting that for $x\in \BR^d\setminus A_J,$
\[
g_J(x)\leq \tilde{\Psi}_{A_o}(o)\leq C I_\alpha \alpha^{-d}
\]
so that by Lemma \ref{lemma:bndedfield}, $g_J(x)$ is uniformly 
bounded by $CI_\alpha/\alpha^d$ where $CI_\alpha<\infty$ is as in 
that lemma. We can therefore use that $e^x-1\leq 2x$ 
for $x\leq 1,$ to see that for $0<s\leq \alpha^d/(CI_\alpha),$ we have
$e^{sg_J(x)}-1\leq 2sg_J(x).$
Hence,
\begin{eqnarray} \label{eqn:gboundeta}
\lefteqn{\int_{\BR^d\setminus A_J}e^{sg_J(x)}-1 dx 
\leq \int_{\BR^d\setminus A_J}2sg_J(x) dx}\\
& & = 2s\int_{\BR^d\setminus A_J}\sum_{j\in J} 
\tilde{l}_\alpha(|y_{c}(B_j)-x|) dx
\leq 2s\sum_{j\in J}\int_{\BR^d\setminus A(B_j)}
\tilde{l}_\alpha(|y_{c}(B_j)-x|) dx \nonumber 
=2sD|J|, \nonumber
\end{eqnarray}
where by using that the side length of $A(B_j)$ is 
$\alpha(4\lceil\sqrt{d}\rceil+1)$, we have that 
\begin{eqnarray*}
\lefteqn{D=\int_{\BR^d\setminus A(B_j)}\tilde{l}_\alpha(|y_{c}(B_j)-x|) dx}\\
& & =\int_{\BR^d\setminus A_o}\tilde{l}_\alpha(|x|)dx
\leq \int_{\BR^d\setminus A_o}\tilde{l}_\alpha(|x|-\alpha \sqrt{d})dx
\leq CI_\alpha
\end{eqnarray*}
as in the proof of Lemma \ref{lemma:bndedfield}. 
Using \eqref{eqn:gboundeta} in \eqref{eqn:cambeta}, we get that 
\begin{equation} \label{eqn:expetaest}
\BE_{\eta_{A_J^c}}\left[e^{s\sum_{j\in J} 
\tilde{\Psi}_{\eta_{A_J^c}}(y_{c}(B_j))}\right]
\leq e^{2sCI_\alpha|J|}.
\end{equation}

We now turn to the second factor on the right hand side of \eqref{eqn:etaxi}.
Observe that for any $k,$ $J\in \{0,1\}^k$ and $j\in J$ we have that 
\[
\tilde{\Psi}_{\xi_D}(y_{c}(B_j))\leq \tilde{\Psi}_{A_o}(o).
\]
Using Lemma \ref{lemma:bndedfield}, it follows that 
\begin{equation} \label{eqn:expxiest}
\BE_{\xi_D}\left[e^{s\sum_{j\in J} \tilde{\Psi}_{\xi_D}(y_{c}(B_j))}\right]
\leq e^{2sCI_\alpha|J|/\alpha^d}.
\end{equation}
Using \eqref{eqn:expetaest} and \eqref{eqn:expxiest}, with \eqref{eqn:etaxi}
we see from \eqref{eqn:Markov} that 
\[
\BP(\tilde{\Psi}(y_c(B_1))\geq h,\ldots,\tilde{\Psi}(y_c(B_k))\geq h|\CD_J)
\leq 
e^{-h|J|s}e^{2sCI_\alpha|J|}e^{2s|J|CI_\alpha /\alpha^d}.
\]
Inserting this into \eqref{eqn:prelest}
\begin{eqnarray*} 
\lefteqn{\BP(\sup_{y\in B_1}\Psi(y)\geq h, \ldots, \sup_{y\in B_k}\Psi(y)\geq h)}\\
& & \leq \sum_{|J|\geq k/2}
e^{-h|J|s}e^{2sCI_\alpha|J|}e^{2s|J|CI_\alpha/\alpha^d}
\BP(\CD_J)+e^{-d(\alpha)k} \nonumber \\
& & \leq e^{-C'I_\alpha hsk/\alpha^d}+e^{-d(\alpha)k},
\end{eqnarray*}
for some $C'>0.$ Finally, by first letting $\alpha$ be so small 
that $e^{-d(\alpha)}\leq \epsilon/2,$
and then taking $h$ large enough, \eqref{eqn:epsest} follows.
\fbox{}\\

We will now prove Theorem \ref{thm:hcfinite} in its entirety. \\

\medskip
{\bf Proof of Theorem \ref{thm:hcfinite}.}
The first statement is simply Theorem \ref{thm:hcfinite_aux} and so we 
turn to the second statement.

Consider the auxiliary attenuation function $l'(r):=l(r+1)$, and let
$\Psi'$ denote the corresponding random field. We observe that 
for any $y\in B(o,1)$ and $x\in \BR^d,$ 
$l'(|x|)=l(|x|+1)\leq l(|x-y|-|y|+1)\leq l(|x-y|),$
so that 
\[
\Psi'(o)=\sum_{x\in \eta}l'(|x|)
\leq \inf_{y\in B(o,1)}\sum_{x\in \eta}l(|x-y|)
=\inf_{y\in B(o,1)}\Psi(y).
\] 
We proceed to show that $\BP(\Psi'(o)=\infty)=1,$ since then it follows that 
$\BP\left(\inf_{y\in \BR^d}\Psi(y)=\infty\right)=1$ by a standard 
countability argument. Therefore, let $A_0:=B(o,1),$ and 
$A_k:=B(o,k+1)\setminus B(o,k)$ and note that 
$Vol(A_k)=\kappa_d((k+1)^d-k^d)\geq d\kappa_d k^{d-1},$
where $\kappa_d$ denotes the volume of the unit ball in dimension $d.$
Furthermore, let $\CA_k$ denote the event that $\eta(A_k)\geq \kappa_d k^{d-1}.$
For any $X\sim Poi(\lambda),$ a standard Chernoff type bound yields
\[
\BP(X\leq \lambda/2)\leq \frac{e^{-\lambda}\left(e\lambda\right)^{\lambda/2}}
{\left(\lambda/2\right)^{\lambda/2}}
=\left(\frac{e}{2}\right)^{-\lambda/2}=e^{-c\lambda},
\]
for some $c>0.$ Therefore, $\BP(\CA_k^c)\leq e^{-cd\kappa_d k^{d-1}}$ so that 
$\BP(\CA_k^c \textrm{ i.o.})=0$ by the Borell-Cantelli lemma.
Thus, for a.e. $\eta,$ there exists a $K=K(\eta)<\infty,$ so that 
$\CA_k$ occurs for every $k\geq K.$ Furthermore, we have that 
if $\CA_k$ occurs, then for any $k\geq 3,$
\begin{eqnarray*}
\lefteqn{\sum_{x\in \eta(A_k)}l'(|x|)\geq \kappa_d k^{d-1}l'(k+1)}\\
& & =\kappa_d k^{d-1}l(k+2)
\geq \kappa_d \frac{k^{d-1}}{(k+3)^{d-1}}\int_{k+2}^{k+3}r^{d-1}l(r)dr
\geq \frac{\kappa_d}{2}\int_{k+2}^{k+3}r^{d-1}l(r)dr.
\end{eqnarray*}
Therefore we get that by letting $K\geq 3,$
\[
\Psi'(o)=\sum_{x\in \eta}l'(|x|)
\geq \sum_{k=K}^\infty \frac{\kappa_d}{2}\int_{k+2}^{k+3}r^{d-1}l(r)dr
=\frac{\kappa_d}{2}\int_{K+2}^{\infty}r^{d-1}l(r)dr=\infty.
\]
\fbox{}\\

\section{Continuity of the field $\Psi$} \label{sec:cont}

In this section we will prove Proposition
\ref{prop:contfield}.
We will often use the following well known equality (see for instance \cite{Kingman} p. 28)
\begin{equation} \label{eqn:VanCamp}
\BE\left[\sum_{x\in \eta} g(x)\right]=\int_{\BR^d} g(x)\mu(dx),
\end{equation}
where $\eta$ is a Poisson process in $\BR^d$ with intensity measure $\mu.$ \\

\noindent
{\bf Proof of Proposition \ref{prop:contfield}.}
We start by proving the statement in the case when $l$ is bounded.
Fix $\alpha,\epsilon>0$, let $g_{y,z}(x)=|l(|x-y|)-l(|x-z|)|,$  
and let $\{D_n\}_{n\geq 1}$ be a sequence of bounded subsets of $\BR^d$
such that $D_n \uparrow \BR^d.$
Observe that for any $\delta>0,$
\begin{eqnarray} \label{eqn:alphaprel}
\lefteqn{\BP\left(\sup_{y,z\in B(o,1):|y-z|<\delta} |\Psi(y)-\Psi(z)|\geq \epsilon\right)}\\
& & \leq 
\BP\left(\sup_{y,z\in B(o,1):|y-z|<\delta} \sum_{x\in \eta} g_{y,z}(x)\geq \epsilon\right) \nonumber\\
& & \leq 
\BP\left(\sup_{y,z\in B(o,1):|y-z|<\delta} \sum_{x\in \eta(D_n)} g_{y,z}(x)\geq \epsilon/2\right)
+\BP\left(\sup_{y,z\in B(o,1):|y-z|<\delta} \sum_{x\in \eta(D_n^c)} g_{y,z}(x)\geq \epsilon/2\right).\nonumber
\end{eqnarray}
We will proceed by bounding the two terms on the right hand side of \eqref{eqn:alphaprel}.
Consider therefore the second term 
\begin{eqnarray} \label{eqn:alpha1}
\lefteqn{
\BP\left(\sup_{y,z\in B(o,1):|y-z|<\delta} \sum_{x\in \eta(D_n^c)} g_{y,z}(x)\geq \epsilon/2\right)}\\
& & \leq \BP\left(\sup_{y,z\in B(o,1):|y-z|<\delta} \sum_{x\in \eta(D_n^c)} l(|x-y|)+l(|x-z|)\geq \epsilon/2\right) \nonumber\\
& & \leq \BP\left(\sup_{y\in B(o,1)} \sum_{x\in \eta(D_n^c)} l(|x-y|)+
\sup_{z\in B(o,1)} \sum_{x\in \eta(D_n^c)}l(|x-z|)\geq \epsilon/2\right) \nonumber\\
& & =\BP\left(\sup_{y\in B(o,1)} \sum_{x\in \eta(D_n^c)} l(|x-y|)\geq \epsilon/4\right).
\nonumber
\end{eqnarray}
Furthermore, we have that for any 
$\epsilon>0,$
\begin{eqnarray} \label{eqn:Dnlimit}
\lefteqn{\BP\left(\sup_{y\in B(o,1)} \sum_{x\in \eta(D_n^c)}l(|x-y|)\geq \epsilon\right)
\leq \frac{1}{\epsilon}\BE\left[\sup_{y\in B(o,1)} \sum_{x\in \eta(D_n^c)}l(|x-y|)\right]}\\
& & \leq \frac{1}{\epsilon}\BE\left[ \sum_{x\in \eta(D_n^c)}\sup_{y\in B(o,1)}l(|x-y|)\right]
=\frac{1}{\epsilon} \int_{\BR^d\setminus D_n} \sup_{y\in B(o,1)}l(|x-y|) dx \nonumber \\
& & \leq \frac{1}{\epsilon} \int_{\BR^d\setminus D_n} l\left(\max(|x|-2,0)\right) dx
\to 0 \textrm{ as } n\to \infty, \nonumber
\end{eqnarray}
where we use \eqref{eqn:VanCamp} in the equality and the fact 
that the intensity measure of $\eta(D_n^c)$ is Lebesgue measure outside
of $D_n.$ We also use the integrability assumption 
$\int_0^\infty r^{d-1}l(r)dr<\infty$ when taking the limit.
By combining \eqref{eqn:alpha1} and \eqref{eqn:Dnlimit}, we see that by taking 
$n$ large enough, the second term of \eqref{eqn:alphaprel} is smaller than $\alpha.$

For the first term, we get that 
\begin{eqnarray} \label{eqn:alpha2}
\lefteqn{
\BP\left(\sup_{y,z\in B(o,1):|y-z|<\delta} \sum_{x\in \eta(D_n)} g_{y,z}(x)
\geq \epsilon/2\right)}\\
& & \leq \frac{1}{\epsilon}
\BE\left[ \sum_{x\in \eta(D_n)}\sup_{y,z\in B(o,1):|y-z|<\delta} |l(|x-y|)-l(|x-z|)|\right]
\nonumber \\
& & =\frac{1}{\epsilon}
\int_{D_n}\sup_{y,z\in B(o,1):|y-z|<\delta} |l(|x-y|)-l(|x-z|)| dx.
\nonumber
\end{eqnarray}
Since $D_n$ is bounded, we have that for any $x\in D_n,$
\[
\sup_{y,z\in B(o,1):|y-z|<\delta} |l(|x-y|)-l(|x-z|)|
\leq \sup_{(r_1,r_2)\in E_n}(l(r_1)-l(r_2))
\]
where $E_n=\{(r_1,r_2)\in \BR^2:0\leq r_1<r_2\leq 2{\rm diam}(D_n),|r_1-r_2|<\delta\}$. Since 
$l(r)$ is continuous, it is uniformly continuous on $[0,2{\rm diam}(D_n)]$.
Therefore, for any fixed $n,$ the right hand side of \eqref{eqn:alpha2} is 
smaller than $\alpha$ for $\delta$ small enough.

We conclude that for any $\epsilon,\alpha>0,$ there exists 
$\delta>0,$ small enough so that
\begin{equation} \label{eqn:2alpha}
\BP\left(\sup_{y,z\in B(o,1):|y-z|<\delta} |\Psi(y)-\Psi(z)|\geq \epsilon\right)
\leq 2 \alpha.
\end{equation}
To conclude the proof, assume that $\Psi(y)$ is not a.s. 
continuous everywhere. Then, with positive probability, there
exists $\epsilon>0$ and a point $w\in B(o,1/2)$ such that for 
any $\delta>0$
\[
\sup_{y:|y-w|<\delta} |\Psi(y)-\Psi(w)|\geq \epsilon,
\]
contradicting \eqref{eqn:2alpha}.

We now turn to the case where $l$ is unbounded.
Then, for any $M<\infty,$ we let $l_M(r)=\min(l(r),M),$
and define $\Psi_M(y)$ to be the random field obtained by using $l_M$ instead of $l.$
If we let  
\[
B_M(x)=\{y\in \BR^d:l(|x-y|)\geq M\},
\]
we see that $\Psi_M(y)=\Psi(y)$ for every $y\in \BR^d\setminus\cup_{x\in \eta}B_M(x).$
By the first case, $\Psi_M(y)$ is continuous everywhere,
and so  $\Psi(y)$ is continuous for any 
$y\in \BR^d\setminus\cup_{x\in \eta}B_M(x).$
Since $M<\infty$ was arbitrary, the statement follows after observing that 
$\lim_{y \to x}\Psi(y)=\infty$ whenever $x\in \eta.$
\fbox{}\\

\section{Uniqueness} \label{sec:uniqueness}

In this section we restrict ourselves to $d=2$. We will first 
consider the case of $\Psi_{> h}$, and then explain how the second
case of Theorem \ref{thm:uniqueness} quickly follows from it. For convenience
we formulate the following separate statement.

\begin{theorem}
\label{groter}
Let $h$ be such that $\Psi_{> h}$ contains an unbounded component.
If $l$ is continuous and with unbounded support, then for $d=2$,
there is a unique such unbounded component.
\end{theorem} 
Our strategy will be to adapt the argument in \cite{GKR} which proves uniqueness 
for a class of models on $\BZ^d$. In order to perform this adaptation
it is much easier to work with {\em arcwise} connectedness rather than 
connectedness. The reason for this is that we can easily form new arcs from intersecting 
arcs, while the corresponding result for connectedness is rather challenging topologically.

However, in our continuous context, we have defined percolation in terms of 
connectedness, as is usually done. But, since $\Psi$ is a.s.\ continuous by Proposition 
\ref{prop:contfield}, the set $\psgh$ is a.s.\ an open set. For 
open sets, connected and arcwise connected are the same thing, as is 
well known. Hence, if $x$ and $y$ are in the same connected component 
of $\psgh$, then there is a continuous curve from $x$ to $y$ in 
$\psgh$. This observation makes $\psgh$ easier to study than $\psgeh$ 
directly, and is the reason for proving Theorem \ref{groter} separately.

In the sequel we try to balance between 
the fact that we do not want or need to repeat the whole argument of 
\cite{GKR} on the 
one hand, and the need to explain in detail what changes are to be made and 
what these changes constitute on the other hand.

In \cite{GKR}, uniqueness of the infinite cluster in two-dimensional 
discrete site percolation is proved under four conditions. Consider a probability 
measure $\mu$ on $\{0,1\}^{\mathbb{Z}^2}$ and let 
$\omega\in \{0,1\}^{\mathbb{Z}^2}$ be a configuration. If $\omega(x)=1$
we call $x\in \BZ^2$ {\em open} and if $\omega(x)=0$ we call it {\em closed}.
The four conditions which together imply uniqueness of the infinite 
open component are:

\begin{enumerate}
\item $\mu$ is invariant under horizontal and vertical translations and under axis reflections.
\item $\mu$ is ergodic (separately) under horizontal and vertical translations.
\item $\mu(E \cap F) \geq \mu(E)\mu(F)$ for events $E$ and $F$ which are both increasing or both decreasing (The FKG inequality).
\item The probability that the origin is in an infinite cluster is non-trivial, that is, strictly between 0 and 1.
\end{enumerate}

In our context, conditions analogous to Conditions 1 and 2 clearly hold. Some 
care is needed for Condition 3 though. We will say that an event $E$ is 
increasing if a configuration in $E$ remains in $E$ if we add additional 
points to the point process $\eta$ (and adapt the field accordingly). 
Furthermore, $E$ is decreasing if $E^c$ is increasing. With these 
definitions, one can prove the analogue of the FKG inequality as in 
the proof of Theorem 2.2 in \cite{MR}.

Condition 4 is natural in the discrete context. Indeed, if the probability that the origin is in an infinite cluster is 0, then by translation invariance, no infinite cluster can exist a.s. The case in which the probability that the origin is in an infinite cluster is 1 was excluded only for convenience, and this assumption is not used in the proof in \cite{GKR}. In our continuous context, we need to be slightly more careful. Suppose that $\psgh$ contains an unbounded component with positive probability. Since $\psgh$ is an open set by continuity of the field, any such unbounded component must be open as well. Hence there must be an $\epsilon>0$ and an $x \in \BR^2$ so that $B(x,\epsilon)$ is contained in an infinite component with positive probability, since a countable collection of such balls covers the plane. By translation invariance, the same must then be true for any $x\in \BR^2$. Hence, any point $x \in \BR^2$ is contained in an infinite component with positive probability, and Condition 4 holds. 

Gandolfi, Keane and Russo prove uniqueness by showing that there exists a $\delta>0$ such that any box $B_n=[-n,n]^2$ is surrounded by an open path with probability at least $\delta$. Hence the probability that all such boxes are surrounded by an open path is at least $\delta$, and since the latter event is translation invariant it must have probability one. This ensures uniqueness, as is well known since 1960, see \cite{Harris}. For the proof of Theorem \ref{groter}, we can 
in principle follow the structure of their argument, with a number 
of modifications, as follows.

\medskip
\noindent
{\bf Proof of Theorem \ref{groter}.}
For any set $A\subset \BR^d,$ we will let $\Gamma_A:=\sup_{x\in A} \Psi(x).$

The first step is to show that it is impossible to have percolation in a 
horizontal strip $Q_M$ of the form
$$
Q_M: =\{(x,y) \in \BR^2; -M \leq y \leq M\},
$$
and similarly for vertical strips. In their case this claim simply follows 
from the fact that closed sites exist (by virtue of Condition 4) and 
then it follows from  Condition 3 that there is a positive probability 
that a strip is blocked completely by closed sites. Finally, ergodicity 
(or rather stationarity) shows that a strip is blocked infinitely many 
times in either direction.

Since we work in a continuous setting, this argument does not go through 
immediately. However, we can adapt it to our context. To this end, consider the 
set $C= C(N,M):=[N,N+1] \times [-M, M]$, that is, a vertical ``strip" in $Q_M$.  
Since the field is a.s.\ finite, by deleting points one by one from $\eta$, 
say in increasing order with respect to distance to $C$, we have that 
upon deleting these points, $\Gamma_C \downarrow 0$. Hence, after 
deleting sufficiently many points it must be the case that 
$\Gamma_C <h$, for any given $h >0$. If we let $\CD^o(L)$
denote the event that the contribution of points outside the 
box $B_{L}$ to $\Gamma_C$ is at most $h$, we conclude that for some
(random) number $L$, $\CD^o(L)$ occurs. It then follows that 
for some deterministic $L_0,$ $\BP(\CD^o(L_0))>0.$ Note also that 
$\CD^o(L_0)$ only depends on the points of $\eta$ outside $B_{L_0}$. 

Let $\CD^i(L_0)$ denote the event that there are no points of $\eta$
in $B_{L_0}$ itself. Then $\BP(\CD^i(L_0))>0$ and by independence of
$\CD^i(L_0)$ and $\CD^o(L_0)$, it also follows that 
$\BP(\CD^i(L_0) \cap \CD^o(L_0))>0$. Furthermore, on the event  
$\CD^i(L_0) \cap \CD^o(L_0),$ we have that $\Gamma_C <h$, and 
we conclude that for any $h>0$ there is positive probability that 
the field $\Psi$ does not exceed $h$ on $C$. 
So any vertical strip in $C(N,M) \subset Q_M$ has positive 
probability to satisfy $\Gamma_C <h$ and by stationarity there must 
be infinitely many of such strips in both directions. So there 
can not be percolation in $Q_M$. 

\medskip

Having established this, we now turn to the second step of the argument.
As mentioned above, in \cite{GKR}, they construct open paths whose union, by virtue of 
their specific construction, surround a given box. They show that 
for any given box, such a construction can be carried out with a 
uniform positive probability. It is at this point of the argument 
that two-dimensionality is crucial as the two-dimensional topology 
forces certain paths to intersect.

The remainder of the proof of uniqueness proceeds in two steps. 
First we prove uniqueness under the assumption that the half-plane 
$H^+:=\{(x,y); y \geq 0\}$ percolates. After that, 
we prove uniqueness under the assumption that $H^+$ does not percolate. 

For $x \in \BR^2, A, B \subset \BR^2$, we write $E(x, A, B)$ for the event 
that there is a continuous path in  $\psgh$ from $x$ to $A$ which is 
contained in $B$, and $E(x, \infty, B)$ if there is an unbounded 
continuous path from $x$ in $B$. We write $L_N:=\{(x,y); y=N\}$,  
$L_N^+:=\{(x,y); y=N, x \geq 0\}$ and  $L_N^-:=\{(x,y); y=N, x \leq 0\}$. 
Finally we write $H_N^+:=\{(x,y); y \geq N\}$, so that $H_0^+ = H^+$.

Let $E :=E(0,\infty, H^+)$, let $D$ be a box centered at the origin, and 
let $D_N := D + (0,N)$. 
Finally, let $\tilde{E}_N := E(0, \infty, H^+ \backslash D_N)$. Now,
\[
\BP(E)\BP(\Gamma_{D}) \geq \BP(\tilde{E}_N)\BP(\Gamma_{D_N}) 
\geq \BP(\tilde{E}_N \cap \Gamma_{D_N})= \BP(E \cap  \Gamma_{D_N}).
\]
Since our system is mixing (see e.g.\ \cite{MR} p. 26
plus the fact that the field is a deterministic function of the Poisson process), 
we have that $\lim_{N \to \infty} \BP(E \cap  \Gamma_{D_N})
=\BP(E)\BP(\Gamma_{D})$. It follows that when $N \to \infty$, $\BP(\tilde{E}_N) \to \BP(E)$. 
In words, if we have percolation from the origin in $H^+$, the conditional probability that 
there is an unbounded path avoiding $D_N$ tends to 1 as $N \to \infty$.

Hence, if the probability that $y_{-N}:=(0,-N)$ percolates in $H^+_{-N}$ is  $\delta >0$, then
for $N$ large enough,
$$
\BP(E(y_{-N}, \infty, H_{-N}^+\backslash D))\geq \delta/2.
$$ 
Since the strip $Q_N$ does not percolate, if $y_{-N}$ percolates in 
$H^+_{-N}\backslash D$, we conclude that the event 
$E(y_{-N}, L_N, Q_N\backslash D)$ must occur, so that 
$\BP(E(y_{-N}, L_N, Q_N\backslash D)) \geq \delta/2$. 

The endpoint of the curve in the event $E(y_{-N}, L_N, Q_N\backslash D)$
is either in $L_N^+$ or in $L_N^-$, and by reflection symmetry, both 
options have the same probability. Hence, 
$$
\BP(E(y_{-N}, L_N^+, Q_N\backslash D))\geq \delta/4.
$$
By reflection symmetry, it then follows that also
$$
\BP(E(y_{N}, L_{-N}^+, Q_N\backslash D))\geq \delta/4,
$$
and by combining the curves in the last two displayed formulas and the 
FKG inequality, we find that
\begin{equation}
\label{curves}
\BP(E(y_{N}, y_{-N}, Q_N\backslash D))\geq \delta^2/16.
\end{equation}
Any curve in the event $E(y_{N}, y_{-N}, Q_N\backslash D)$
either has $D$ on the left or on the right (depending whether it has 
positive or negative winding number) and again by reflection symmetry, 
both possibilities must have probability at least $\delta^2/32$.  
Let $J^+$ ($J^-$) be the sub-event of $E(y_{N}, y_{-N}, Q_N\backslash D)$
where there exists a curve with positive (negative) winding number.
By the FKG inequality, we 
have that $\BP(J^+ \cap J^-) \geq \delta^4/1024$. But on $J^+ \cap J^-$, 
the box $D$ is surrounded by a continuous curve in $\psgh$, and we are done. 

Finally, we consider the case in which the half space does not percolate. 
We can modify the argument in \cite{GKR} similarly and we do not spell out 
all details. In the first case we showed that if we have percolation from 
the origin in $H^+$, the conditional probability that there is an unbounded 
path avoiding $D_N$ tends to 1 as $N \to \infty$. In this second case it 
turns out that we need to show that this remains true if we in addition 
also want to avoid $D_{-N}$. For this, the usual mixing property that we 
used above, does not suffice, and a version of 3-mixing is necessary. As 
in \cite{GKR}, we use Theorem 4.11 in \cite{Fu} for this, in which it is 
shown that ordinary weak mixing implies 3-mixing along a sequence of density 1. 
Since our system is weakly mixing, this application of Theorem 4.11 in 
\cite{Fu} is somewhat simpler than in \cite{GKR}, but other than that 
our argument is the same, and we do not repeat it here.
\fbox{}\\

\medskip
Finally we show how Theorem \ref{groter} implies Theorem \ref{thm:uniqueness}.

\medskip
\noindent
{\bf Proof of Theorem \ref{thm:uniqueness}}.
We first claim that $\Psi_{>h}$ percolates if and only if $\Psi_{\geq h}$ percolates. The ``only if" is 
clear, since $\Psi_{> h} \subset \Psi_{\geq h}$. 

Next, suppose that $\Psi_{\geq h}$ percolates. By definition, this implies that $\psgeh$ a.s.\ contains an unbounded connected component. Let us denote this event by $C_h$. Let $A$ be a bounded region with positive volume. Since the probability of $C_h$ is 1, it must be the case that 
$$
\BP(C_h| \eta(A)=0)=1,
$$
where $\eta(A)$ is the number of points of the Poisson process in $A$.
Since we can sample from the conditional distribution of the process given $\eta(A)=0$ by first sampling unconditionally and then simply remove all points in $A$, it follows that we cannot destroy the event of percolation by removing all points in $A$.

Hence if we take all points out from $A$, the resulting field $\Psi^A$, say, will be such that $\Psi^A_{\geq h}$ percolates a.s. But if $\eta(A) >0$, then it is the case that 
$\Psi^A_{\geq h} \subseteq \psgh$, and it is precisely here we assume that 
the attenuation function $l$ has unbounded support. Hence, with positive probability we have that $\psgh$ percolates, and by ergodicity this implies that $\psgh$ contains an infinite component a.s.

We can now quickly finish the proof. Suppose that $\psgeh$ percolates. Then, as we just saw, also 
$\psgh$ percolates. Hence we can apply the proof of Theorem \ref{groter}, and conclude that 
$\psgh$ does contain continuous curves around each box. Since $\psgeh$ is an even larger set, the 
same must be true for $\psgeh$ and uniqueness for this latter set follows as before.
\fbox{}\\

\noindent 
{\bf Remark:}In light of Theorem \ref{thm:uniqueness}, it is of course natural to expect 
that uniqueness should hold also for $d\geq 3.$ The classical argument for uniqueness
in various lattice models and continuum percolation consists of two parts. Below
we examine these separately.


Let $N_h$ be the number of unbounded components in $\Psi_{\geq h}$. 
Following the arguments of \cite{NS} (which is for the lattice case 
but can easily be adapted to the setting of Boolean percolation, see 
\cite{MR} Proposition 3.3) one starts by observing that
$\BP(N_h=k)=1$ for some $k\in\{0,1,\ldots\}\cup\{\infty\}.$  
Assume for instance that $\BP(N_h=3)=1$, and proceed by taking a 
box $[-n,n]^d$ large enough so that 
at least two of these
infinite components intersect the box with positive probability. Then, glue these 
two components together by the use of a finite energy argument. That is, turn all 
sites in the box to state $1$ in the discrete case, or add balls to 
the box in the 
Boolean percolation case. In this way, we reduce $N_h$ by (at least) 1, 
showing that 
$\BP(N_h=3)<1,$ a contradiction. If one attempts to repeat this procedure 
in our setting (with the support of $l$ being unbounded), 
one finds that by adding points to the field, the gluing of two 
infinite components might at the same
time result in the forming of a completely new infinite component 
somewhere outside the box. 
Therefore, one cannot conclude that $\BP(N_h=3)<1$.

The second difficulty occurs when attempting to rule out the possibility that 
$\BP(N_h=\infty)=1.$ For the Boolean percolation model one uses an argument by Burton and Keane in \cite{BK},
adapted to the case of Boolean percolation (see \cite{MR} proof of Theorem 3.6). 
However this argument hinges on the trivial
but crucial fact that for this model any unbounded component must contain infinitely many 
points of the Poisson process $\eta$. This is not the case in our setting. An unbounded component
can in principle contain only a finite number of points of $\eta$, 
or indeed none at all. \\

We now turn to the last result of this paper, Theorem \ref{thm:thetacont}.
In order to prove continuity, we will give separate arguments for 
left- and right-continuity. The strategy to prove right-continuity 
will be similar to the corresponding result (i.e. left-continuity)
for discrete lattice percolation (see \cite{Grimmett}, Section 8.3).
However, while the other case is trivial for discrete percolation, 
this is where most of the effort in proving Theorem \ref{thm:thetacont}
lies. Before giving the full proof, we will need to establish two 
lemmas that will 
be used to prove left-continuity. See also the remark after the end of 
the proof of Theorem \ref{thm:thetacont}.

Let $X^0\sim$Poi$(\lambda)$, and let $X^1=X^0+1.$ The following
lemma provides a useful coupling.
\begin{lemma} \label{lemma:Poicouple}
There exist random variables $Y^0\stackrel{d}{=}X^0$ and 
$Y^1\stackrel{d}{=}X^1$ coupled so that 
\[
\BP(Y^0\neq Y^1)
=\frac{\lambda^{\lfloor \lambda \rfloor}+1}{(\lfloor \lambda \rfloor+1)!}
e^{-\lambda}.
\]
\end{lemma}
\noindent
{\bf Proof.}
In what follows, sums of the form $\sum_{l=M}^{M-1} a_l$ is understood to be 0,
and in order not to introduce cumbersome notation, expressions such as
$\lambda^k/k!$ will be interpreted as 0 for $k<0.$
Note also that 
\begin{equation} \label{eqn:Poiratio}
\frac{\BP(X^0=k)}{\BP(X^1=k)}
=\frac{\frac{\lambda^k}{k!}e^{-\lambda}}
{\frac{\lambda^{k-1}}{(k-1)!}e^{-\lambda}}
=\frac{\lambda}{k}\geq 1 \textrm{ iff } k\leq \lambda.
\end{equation}

We start by giving the coupling and then verify that it is well defined
and has the correct properties.
Let $U \sim U[0,1]$ and for $1\leq k \leq \lambda$ let $Y^0=Y^1=k$ if
\begin{equation} \label{eqn:Poi1}
\sum_{l=0}^{k-2}\frac{\lambda^{l}}{l!}e^{-\lambda}
<U\leq \sum_{l=0}^{k-1}\frac{\lambda^{l}}{l!}e^{-\lambda}
\end{equation}
while for $k>\lambda$ we let $Y^0=Y^1=k$ if
\begin{equation} \label{eqn:Poi2}
\sum_{l=0}^{\lfloor\lambda\rfloor}\frac{\lambda^l}{l!}e^{-\lambda}
+\sum_{l=\lfloor\lambda\rfloor +2}^{k-1}\frac{\lambda^l}{l!}e^{-\lambda}
<U\leq \sum_{l=0}^{\lfloor\lambda\rfloor}\frac{\lambda^l}{l!}e^{-\lambda}
+\sum_{l=\lfloor\lambda\rfloor +2}^{k}\frac{\lambda^l}{l!}e^{-\lambda}
\end{equation}
Furthermore, we let $Y^0=k$ if $k\leq \lambda$ and
\begin{equation} \label{eqn:Poi3}
1-\frac{\lambda^k}{k!}e^{-\lambda}
< U\leq 1-\frac{\lambda^{k-1}}{(k-1)!}e^{-\lambda}
\end{equation}
while $Y^1=k$ if $k> \lambda$ and
\begin{equation} \label{eqn:Poi4}
1-\frac{\lambda^{k-1}}{(k-1)!}e^{-\lambda}<U
\leq 1-\frac{\lambda^{k}}{k!}e^{-\lambda}.
\end{equation}
Consider now \eqref{eqn:Poi3}. Since $k\leq \lambda$ it follows from 
\eqref{eqn:Poiratio} that 
\[
1-\frac{\lambda^k}{k!}e^{-\lambda} 
\leq 1-\frac{\lambda^{k-1}}{(k-1)!}e^{-\lambda},
\]
with equality iff $k=\lambda.$ It follows that \eqref{eqn:Poi3}
is well defined, and similarly we can verify that \eqref{eqn:Poi4}
is also well defined.

We proceed to verify that \eqref{eqn:Poi1}--\eqref{eqn:Poi4} 
gives the correct distributions of $Y^0$ and $Y^1$. 
To this end, observe that from \eqref{eqn:Poi1} and \eqref{eqn:Poi3} 
we have that for $0\leq k\leq \lambda$ 
\[
\BP(Y^0=k)=\frac{\lambda^{k-1}}{(k-1)!}e^{-\lambda}
+\frac{\lambda^{k}}{k!}e^{-\lambda}-\frac{\lambda^{k-1}}{(k-1)!}e^{-\lambda}
=\frac{\lambda^{k}}{k!}e^{-\lambda}.
\]
Furthermore, from \eqref{eqn:Poi1} we get that for $k>\lambda,$
\[
\BP(Y^0=k)=\frac{\lambda^{k}}{k!}e^{-\lambda}
\]
so that indeed $Y^0\sim$Poi$(\lambda).$

Similarly, we see from \eqref{eqn:Poi1} that for 
$1\leq k \leq \lambda$ we have that 
\[
\BP(Y^1=k)=\frac{\lambda^{k-1}}{(k-1)!}e^{-\lambda},
\]
while by summing the contributions from \eqref{eqn:Poi2}¨and 
\eqref{eqn:Poi4} we get that for $k>\lambda$
\[
\BP(Y^1=k)
=\frac{\lambda^k}{k!}e^{-\lambda}+\frac{\lambda^{k-1}}{(k-1)!}e^{-\lambda}
-\frac{\lambda^{k}}{k!}e^{-\lambda}
=\frac{\lambda^{k-1}}{(k-1)!}e^{-\lambda},
\]
so that $Y^1$ has the desired distribution. 

Finally, the lemma follows by observing that 
\[
\BP(Y^0\neq Y^1)=\BP\left(U>\sum_{l=0}^{\lfloor\lambda\rfloor}\frac{\lambda^l}{l!}e^{-\lambda}
+\sum_{l=\lfloor\lambda\rfloor+2}^{\infty}\frac{\lambda^l}{l!}e^{-\lambda}\right)
=\frac{\lambda^{\lfloor \lambda \rfloor}+1}{(\lfloor \lambda \rfloor+1)!}
e^{-\lambda}.
\]
\fbox{}\\

Let $\eta^0_n$ be a homogeneous Poisson process in $\BR^2$ with rate 
1, and let $\eta^1_n$ be a point process such that 
$\eta^1_n\stackrel{d}{=} \eta^0_n+\delta_{V_n}$ where $V_n\sim$U$(B_n)$
and $B_n=[-n,n]^2$.
Thus $\eta^1_n$ is constructed by adding a point uniformly located within 
the box $B_n$ to a homogeneous Poisson process in $\BR^2$.
Let $\BP_n^i$ be the distribution of $\eta^i_n$ for $i=0,1$
and let 
\[
d_{TV}(\BP_n^0,\BP_n^1):=\sup_{A}|\BP_n^0(A)-\BP_n^1(A)|
\]
be the total variation distance between $\BP_n^0$ and $\BP_n^1$, where the supremum
is taken over all measurable events $A$.

\begin{lemma} \label{lemma:PoiTV}
For any $n\geq 1$ we have that 
\[
d_{TV}(\BP_n^0,\BP_n^1)\leq \frac{(4n^2)^{4n^2+1}}{(4n^2+1)!}e^{-4n^2}
\leq n^{-1}.
\]
\end{lemma}
\noindent
{\bf Proof.}
Let $\lambda=4n^2$ and pick $Y^0,Y^1$ as in Lemma \ref{lemma:Poicouple}.
Furthermore, let $\eta$ be a homogeneous Poisson process in $\BR^2,$
independent of $Y^0$ and $Y^1,$
and let $(U_k)_{k\geq 1}$ be an i.i.d. sequence independent of $Y^0,Y^1$
and $\eta$ and such that $U_k\sim$U$(B_n)$. Then, define 
\[
\eta_n^0:=\eta(B_n^c)+\sum_{k=1}^{Y^0}\delta_{U_k},
\]
and
\[
\eta_n^1:=\eta(B_n^c)+\sum_{k=1}^{Y^1}\delta_{U_k}.
\]
It is easy to see that $\eta_n^i\sim \BP_n^i$ and that 
\begin{equation} \label{eqn:etaY}
\BP(\eta_n^0\neq \eta_n^1)=\BP(Y^0\neq Y^1).
\end{equation}
Thus, for any measurable event $A$, we have that 
\[
|\BP_n^0(A)-\BP_n^1(A)|=\BP(\eta_n^0\in A, \eta_n^1 \not \in A)
+\BP(\eta_n^1\in A, \eta_n^0 \not \in A)\leq 
\BP(\eta_n^0\neq \eta_n^1)
\leq \frac{(4n^2)^{4n^2+1}}{(4n^2+1)!}e^{-4n^2},
\]
by using \eqref{eqn:etaY} and Lemma \ref{lemma:Poicouple}.

Furthermore, by Stirling's approximation, we see that 
\[
\frac{(4n^2)^{4n^2+1}}{(4n^2+1)!}e^{-4n^2}
\leq \frac{(4n^2)^{4n^2}}{4n^2!}e^{-4n^2}
\leq \frac{(4n^2)^{4n^2}}{\sqrt{2\pi}(4n^2)^{4n^2+1/2}e^{-4n^2}}e^{-4n^2}
\leq n^{-1}.
\]
\fbox{}\\
\noindent
{\bf Remark:} Although we choose to state and prove this only for $d=2,$
a version of this lemma obviously holds for all $d\geq 1.$ \\

We are now ready to give the proof of our last result.\\
\noindent
{\bf Proof of Theorem \ref{thm:thetacont}.} 
We start by proving the left-continuity of $\theta_{>}(h).$ 
We claim that 
\begin{equation} \label{eqn:uplimg1}
\lim_{g \uparrow h} \theta_{>}(g)=\BP(\CC_{o,>}(g) \textrm{ is unbounded for every }
g<h)=\BP(\CC_{o,\geq}(h) \textrm{ is unbounded}).
\end{equation}
To see this, observe first that trivially 
\[
\{\CC_{o,\geq}(h) \textrm{ is unbounded}\} \subset 
\bigcap_{g<h} \{\CC_{o,>}(g) \textrm{ is unbounded}\}.
\]
Secondly, assume that $\CC_{o,\geq}(h)$ is bounded. Since 
$\CC_{o,\geq}(h)$ and $\Psi_{\geq h} \setminus \CC_{o,\geq}(h)$
are disconnected, there exist open sets $G_1,G_2$ such that $G_1$ is connected,
$\CC_{o,\geq}(h) \subset G_1$, $\Psi_{\geq h} \setminus \CC_{o,\geq}(h)\subset G_2$
and $G_1 \cap G_2=\emptyset.$ Therefore, the set $G_3=G_1\setminus \CC_{o,\geq}(h)$
is an open connected set separating the origin $o$ from $\infty.$ Since 
$G_3$ is then also arcwise connected, it follows that it 
 must contain a circuit surrounding the origin.
That is, there exists a continuous function $\gamma:[0,1] \to \BR^2$ such 
that $\gamma(0)=\gamma(1)$ and $\gamma$ separates $o$ from $\infty.$
Since $\gamma$ is continuous,
the image of $\gamma$ (Im$(\gamma)$) is compact, and so 
$\sup_{t\in[0,1]}\Psi(\gamma(t))$ is obtained, since $\Psi$ is continuous
by Proposition \ref{prop:contfield}. By construction, $G_3 \subset \BR^2\setminus 
\Psi_{\geq h}$ and so 
Im$(\gamma) \subset \BR^2 \setminus \Psi_{\geq h}.$ We conclude that 
$\sup_{t\in[0,1]}\Psi(\gamma(t))<h.$ Therefore, for any $g$ such that
$\sup_{t\in[0,1]}\Psi(\gamma(t))<g<h$ we must have that $\CC_{o,>}(g)$ is bounded.
This proves \eqref{eqn:uplimg1}.

Let $n$ be any integer and take
\[
\eta\in \{\CC_{o,\geq}(h) \textrm{ is unbounded}\} \setminus 
\{\CC_{o,>}(h) \textrm{ is unbounded}\}.
\]
Let $\eta^1_n=\eta+\delta_{V_n}$ where $V_n\sim$U$(B_n)$ and observe 
that since $l$ has unbounded support,
\[
\eta_n^1\in \{\CC_{o,>}(h) \textrm{ is unbounded}\}.
\]
Using Lemma \ref{lemma:PoiTV} we get that 
\begin{eqnarray*}
\lefteqn{\BP(\CC_{o,\geq}(h) \textrm{ is unbounded})} \\
& & \leq \BP^1_n(\CC_{o,>}(h) \textrm{ is unbounded})
\leq \BP^0_n(\CC_{o,>}(h) \textrm{ is unbounded})+n^{-1}
=\theta_{>}(h)+n^{-1}.
\end{eqnarray*}
This together with \eqref{eqn:uplimg1} yields
\[
\lim_{g \uparrow h} \theta_{>}(g)\leq \lim_{n \to \infty}\theta_{>}(h)+n^{-1}
=\theta_{>}(h).
\]

It remains to prove that 
\begin{equation} \label{eqn:downlimg1}
\lim_{g \downarrow h} \theta_{>}(g)=\theta_{>}(h),
\end{equation}
for $h<h_c,$ and we will use a similar approach as above.
 We note that 
\begin{equation} \label{eqn:downlimg2}
\lim_{g \downarrow h} \BP(\CC_{o,>}(g) \textrm{ is unbounded})
=\BP(\CC_{o,>}(g) \textrm{ is unbounded for some } g>h).
\end{equation}
Assume that $\CC_{o,>}(h)$ contains an unbounded component, and consider any 
$h<g<h_c.$ Since $\Psi_{>g}$ also must contain an unbounded component $I_g,$
and since by Theorem \ref{thm:uniqueness} we know that this is unique, 
we conclude that $I_g \subset \CC_{o,>}(h).$ As above, $\CC_{o,>}(h)$
is an open set, and therefore arcwise connected. Thus, for $z\in I_g,$
there exists a continuous function $\phi:[0,1]\to \BR^2$ 
such that $\phi(0)=o$ and $\phi(1)=z.$ Since $\phi$ is continuous,
the Im$(\phi)$ is compact, and so 
$\inf_{t\in[0,1]}\Psi(\phi(t))$ is obtained, since $\Psi$ is continuous
by Proposition \ref{prop:contfield}. Furthermore, since 
Im$(\phi) \subset \CC_{o,>}(h),$ we conclude that 
$\inf_{t\in[0,1]}\Psi(\phi(t))>h.$ Therefore, for some 
$h<g'<h_c$ we also have that $\inf_{t\in[0,1]}\Psi(\phi(t))>g'$,
and so $\CC_{o,>}(g')$ contains an unbounded component.
We conclude that 
\begin{equation} \label{eqn:someg}
\BP(\CC_{o,>}(g) \textrm{ is unbounded for some } g>h)
=\BP(\CC_{o,>}(h) \textrm{ is unbounded}).
\end{equation}
Combining equations \eqref{eqn:downlimg2} and \eqref{eqn:someg} we conclude
that \eqref{eqn:downlimg1} holds.

In order to complete the proof, we simply observe that for any 
$g<h$ we have that $\theta_>(h)\leq \theta_{\geq}(h)\leq \theta_>(g)$
so that 
\[
\theta_>(h)\leq \theta_{\geq}(h) 
\leq \liminf_{g \uparrow h} \theta_{>}(g)=\theta_>(h),
\]
so that indeed $\theta_>(h)=\theta_{\geq}(h)$ for every $h<h_c.$
\fbox{}\\

\noindent
{\bf Remark:} 
%
Consider the event $\{o\leftrightarrow \partial B_n\}$, that the 
origin is connected to the boundary of $B_n.$
In the discrete case, it is trivial that 
$\BP_p(o \leftrightarrow \partial B_n)$ is continuous as a function of the percolation
parameter $p,$  since it is an event that depends on the state of only finitely many
edges. This then gives an easy proof of right-continuity (corresponding to left-continuity 
in our case). In our model, points of $\eta$ at any distance contribute to 
the field in $B_n.$ Therefore, we cannot claim immediate continuity of 
$\BP(o \leftrightarrow \partial B_n)$, although our methods above can be used to prove
it.

\end{document}